\documentclass[11pt]{article}

\def\init{\setcounter{equation}{0}}
\newtheorem{theorem}{Theorem}[section]

\newtheorem{lemma}[theorem]{Lemma}
\newtheorem{definition}[theorem]{Definition}
\newtheorem{corollary}[theorem]{Corollary}
\newtheorem{conjecture}[theorem]{Conjecture}

\newcommand{\qed}{{\vrule height 1ex width 1ex depth -.1ex }}
\def\rr{{\bf R}}
\def\cc{{\bf C}}
\def\nn{{\bf N}}
\def\zz{{\bf Z}}

\def\cala{{\cal A}}
\def\calb{{\cal B}}

\def\jap#1{\langle #1\rangle}

\def\proof{{\it  Proof. }}

\def\eps{\epsilon}

\def\mod{\hbox{mod}\,}

\def\div{\,|\,}
\def\ndiv{\!\not|\,}

\def\eps{\epsilon}

\begin{document}

\title{Spectra of certain types of polynomials and tiling
of integers with translates of finite sets}
\author{Sergei Konyagin and Izabella {\L}aba}
\maketitle


\section{Introduction}
\label{intro}
\init


\begin{definition}\label{N-spectrum}
Let $A(x)\in\zz[x]$ be a polynomial.  We say that $\{\theta_1,
\theta_2,\dots,\theta_{N-1}\}$ is an {\em $N$-spectrum} for $A(x)$ if the
$\theta_j$ are all distinct and
\[
A(\eps_{jk})=0\hbox{ for all }0\leq j,k\leq N-1,\ j\neq k,
\]
where
\[
\eps_{jk}=e^{2\pi i(\theta_j-\theta_k)},\ \theta_0=0.
\]
\end{definition}

Definition \ref{N-spectrum} is motivated by a conjecture of Fuglede \cite{Fug},
which asserts that a measurable set $E\subset\rr^n$ tiles $\rr^n$ by 
translations
if and only if the space $L^2(E)$ has an orthogonal basis consisting of 
exponential functions $\{e^{2\pi i \lambda\cdot x}\}_{\lambda\in\Lambda}$; the set
$\Lambda$ is called a {\em spectrum} for $E$.  For recent work
on Fuglede's conjecture see e.g. \cite{IKP}, \cite{IKT1}, \cite{IKT2}, \cite{IP},
\cite{JP}, \cite{JP2}, \cite{K1}, \cite{K2}, \cite{KP}, \cite{L2}, \cite{L},
\cite{LRW}, \cite{LW2}, \cite{P1}, \cite{P2}, \cite{PW}.  
In the special case when $E\subset\rr$ is a union of $N$ intervals of length 1,
Fuglede's conjecture was proved in \cite{JP} (see also \cite{LW2},
\cite{LS}, \cite{L}, \cite{P1}) to be equivalent to the following.

\begin{conjecture}\label{spectrum-conj}
Let $A(x)$ be a polynomial whose all coefficients are nonnegative integers.
Then the following are equivalent:

\smallskip
(T) There is a finite set $A\subset\{0,1,2,\dots\}$ such that $A(x)=\sum_{a\in A}x^a$.
Moreover, the set
$A$ tiles $\zz$ by translations, i.e., there is a set $B\subset\zz$ 
(called the {\em translation set}) such that every integer $n$ can be 
uniquely represented as $n=a+b$ with $a\in A$ and $b\in B$;

\smallskip
(S) $A(x)$ has an $N$-spectrum with $N=A(1)$.

\end{conjecture}

If (T) holds for some $B$, we will write $A\oplus B=\zz$.
Throughout this paper we will always assume that $2\leq \#A<\infty$.

We will also address the question of characterizing finite sets $A$ which
tile $\zz$ by translations.  This problem has been considered by several
authors and is closely related to many
questions concerning factorization of finite groups, in particular
periodicity and replacement of factors; see e.g. \cite{CM}, \cite{GLW}, \cite{LS},
\cite{LW2}, \cite{New}, \cite{Sands}, \cite{Sands2}, \cite{Tij}.
In particular, the following conditions were formulated in
\cite{CM}.  Let $\Phi_s(x)$ denote the {\em $s$-th cyclotomic polynomial},
defined inductively by
\begin{equation}\label{def-cyclo}
  x^n-1=\prod_{s\div n}\Phi_s(x).  
\end{equation}
We define $S_A$ to be the set of prime powers $p^\alpha$ such that
$\Phi_{p^\alpha}(x)$ divides $A(x)$.

\begin{conjecture}\label{CM-conj} 
$A$ tiles $\zz$ by translations if and only if the following two 
conditions hold:

\medskip
(T1) $A(1)=\prod_{s\in S_A} \Phi_s(1)$,

\medskip
(T2) if $s_1,\dots,s_k\in S_A$ are powers of different primes, then
$\Phi_{s_1\dots s_k}(x)$ divides $A(x)$.
\end{conjecture}

It is proved in \cite{CM} that (T1)-(T2) imply (T), (T) implies (T1),
and that (T) implies (T2) under the additional assumption that $\#A$
has at most two distinct prime factors.  It is not known whether 
(T) always implies (T2); a partial result in the three-prime case
was obtained in \cite{GLW}.  

Conjecture \ref{CM-conj}, if true, implies
one part of Conjecture \ref{spectrum-conj}, since (T1)-(T2) imply (S)
with the spectrum
\[
\Big\{\sum_{s\in S_A}\frac{k_s}{s}:\ 0\leq k_s<p\hbox{ if }s=p^\alpha,\ 
p\hbox{ prime }\Big\}\setminus\{0\}
\]
(see \cite{L}). In particular, we have (T)$\Rightarrow$(S) if $\# A$ has
at most two distinct prime factors.

Conjectures \ref{spectrum-conj} and \ref{CM-conj} have been verified 
in several other special cases.  They are both true under the assumption that 
the degree of $A(x)$ is less than $\frac{3N}{2}-1$, where $N=\#A$ 
\cite{L}.  It also follows from the results of \cite{L2} that Conjecture 
\ref{spectrum-conj} is true for polynomials of the form 
\[
A(x)=\frac{x^{k}-1}{x-1}+x^{m}\frac{x^{n}-1}{x-1}
=1+x+\dots+x^{k-1}+x^{m}+x^{m+1}+\dots+x^{m+n-1}
\]
with $m\geq k$.  It is also known \cite{PW} that if a set $A$
tiles the nonnegative integers by translations, then (S) holds (in fact
the result of \cite{PW} applies to more general sets $E\subset[0,\infty)$).
Finally, it is proved in \cite{L} that if $A(x)$ has degree
less than $\frac{5N}{2}-1$, where $N=\#A$, then an $N$-spectrum must
be rational.

The results of this paper are as follows.

\begin{theorem}\label{irreducible}
Conjectures \ref{spectrum-conj} and \ref{CM-conj} are true if 
$A(x)$ is assumed to be irreducible.
Furthermore, if $A(x)$ is irreducible, then (T),(S) hold if and only if 
$\#A=p$ is prime and $A(x)=1+x^{p^{\alpha-1}}+x^{2p^{\alpha-1}}+\dots
+x^{(p-1)p^{\alpha-1}}$ for some $\alpha\in\nn$.
\end{theorem}

Our next two theorems concerns polynomials of the form
\begin{equation}\label{A1A2}
A(x)=\prod_{i=1}^N A_i(x),\  
A_i(x)=1+x^{m_i}+\dots+x^{m_i(n_i-1)}=\frac{x^{m_in_i}-1}{x^{m_i}-1}.
\end{equation}
Note that each factor $A_i$ is the characteristic polynomial of the set
$\{0,m_i,2m_i,\dots,(n_i-1)m_i\}$, which tiles $\zz$ with the translation
set $\{0,1,\dots,m_i-1\}+ m_in_i\zz$.  Furthermore, $A_i(1)=n_i$ and
each $A_i$ has an $n_i$-spectrum $\{k/n_im_i:\ k=1,2,\dots,n_i-1\}$.
It follows from Corollary \ref{vectors-cor} that $A(x)$ cannot have an 
$M$-spectrum with $M=n_1\dots n_N=A(1)$ unless all coefficients of $A(x)$ are
0 or 1, i.e. $A(x)$ is a characteristic polynomial of a set $A\subset\zz$.

\begin{theorem}\label{two-factors}
Conjectures \ref{spectrum-conj} and \ref{CM-conj} are true for polynomials
of the form (\ref{A1A2}) with $N=2$.
\end{theorem}

\begin{theorem}\label{N-factors}
Conjecture \ref{CM-conj} is true for polynomials of the form (\ref{A1A2})
for all $N\geq 2$.
\end{theorem}

\section{Preliminaries}
\label{prel}
\init

It is well known (see e.g. \cite{New}) that all tilings of $\zz$
by finite sets are periodic: if $A$ is finite and $A\oplus C=\zz$,
then $C=B\oplus M\zz$ for some finite set $B$ such that $\#A\cdot
\#B=M$.  Equivalently, $A\oplus B$ is a complete residue system
modulo $M$, with $M$ as above.  We can rewrite it as
\begin{equation}\label{e.tiling}
  A(x)B(x)=1+x+\dots +x^{M-1}\ (\mod (x^{M}-1)),
\end{equation}
where $B(x)=\sum_{b\in B}x^b$.  By (\ref{def-cyclo}), this is
equivalent to
\begin{equation}\label{e.tiling2}
  A(1)B(1)=M\hbox{ and }\Phi_s(x)\div A(x)B(x)\hbox{ for all }
s\div M,\ s\neq 1.
\end{equation}

The following lemma is due to A. Granville (unpublished).

\begin{lemma}\label{granville}
If $A$ tiles $\zz$ by translations, then it admits a tiling whose
period divides the number
\[L=\hbox{lcm}\{s:\ \Phi_s(x)\div A(x)\}.\]
\end{lemma}

\proof
Fix $A$, and let $A\oplus B=\zz_M\ (\mod M)$.  Replacing $B$
by $\{c\in\{0,1,\dots,M-1\}:\ c=b\ (\mod M)$ for some $b\in B\}$
if necessary, we may assume that $B\subset\{0,\dots,M-1\}$.   
Let $l=(L,M)$.  If $d\div M$ but $d\ndiv L$ then 
\[
\Phi_d(x)\,\Big|\,\frac{x^M-1}{x-1}\,\Big|\, A(x)B(x)
\]
but $\Phi_d(x)\ndiv A(x)$, hence $\Phi_d(x)\div B(x)$.  Therefore
\[
\frac{x^M-1}{x^l-1}=\prod_{d\div M,d\ndiv l}\Phi_d(x)\,\Big|\, B(x).
\]
Let $P(x)=B(x)(x^l-1)/(x^M-1)=\sum_{j=0}^{l-1}p_jx^j$.  Then
$$B(x)=\frac{x^M-1}{x^l-1}P(x)=\sum_{j=0}^{l-1}p_j(x^j+x^{j+l}
+\dots+x^{j+M-l}).$$
Thus the polynomial $P(x)$ has the form $P(x)=B_0(x)$,
where $B_0=\{b\in B:\ 0\leq b\leq l-1\}$.
Then $A(x)B_0(x)=\frac{x^l-1}{x-1}\ (\mod (x^l-1))$ and 
$A(x)B_0(x)=\zz_l\ (\mod l)$.
\qed

\bigskip

We will need the following well known property of cyclotomic polynomials:
\begin{equation}\label{e.cyclo}
  \Phi_s(1)=\left\{
  \begin{array}{ll}
  0\ &\hbox{if }s=1,\\
  p\ &\hbox{if }s=p^\alpha,\ p\hbox{ prime },\\
  1\ &\hbox{otherwise.}\\
  \end{array}\right.
\end{equation} 

Finally, we will need the following lemma.

\begin{lemma}\label{vectors}
Suppose that $A(x)\in\zz[x]$ has nonnegative coefficients.  Then $A(x)$
cannot have an $N$-spectrum for any $N$ greater than the number of
non-zero coefficients of $A$.
\end{lemma}

\proof The proof is a simple modification of an argument of \cite{JP}.
Let $A(x)=\sum_{j=1}^M a_jx^{\alpha_j}$, where $a_j>0$ for all $j$.
Let $\{\theta_j:\ j=1,\dots,N-1\}$ be an $N$-spectrum for $A(x)$,
$\theta_N=0$, $\eps_{j}=e^{2\pi i\theta_j}$ and
$\eps_{jk}=e^{2\pi i(\theta_j-\theta_k)}$. 
Then the condition $A(\eps_{jk})=0$ means that the vectors
\[
{\bf u}_j=(\eps_j^{\alpha_1},\dots,\eps_j^{\alpha_M})
\]
are mutually orthogonal in $\cc^M$ with respect to the inner
product 
\[
({\bf v},{\bf w})=\sum a_kv_kw_k,\ {\bf v}=(v_1,\dots,v_M),\ 
{\bf w}=(w_1,\dots,w_M).
\]
Since there can be at most $M$ such vectors, it follows that
$N\leq M$.
\qed

\begin{corollary}\label{vectors-cor}
Assume that $A(x)\in\zz[x]$ has nonnegative coefficients, and that it
satisfies either (T1)-(T2) or (S).  Then all non-zero coefficients of
$A(x)$ are 1.
\end{corollary}

\proof
If $A(x)$ satisfies (T1)-(T2), then it also satisfies (S) \cite{L}.
Thus it suffices to consider the case when (S) holds.  But then
the corollary is an immediate consequence of Lemma \ref{vectors}.
\qed

\section{Proof of Theorem \ref{irreducible}}
\label{sec2}
\init

Throughout this section we assume that $A(x)$ is irreducible.
Assume that $A$ tiles $\zz$ by translations.  Then (T1) holds, and 
it follows from the irreducibility of $A(x)$ and (\ref{e.cyclo})
that $A(x)=\Phi_{p^\alpha}(x)$
for some prime $p$.  Hence $N=A(1)=p$ and the set 
$\{jp^{-\alpha}:\ j=1,2,\dots,p^\alpha -1\}$ is an $N$-spectrum for $A$.

Suppose now that $A(x)$ has an $N$-spectrum.  Let $e(u)=e^{2\pi iu}$, 
$$A(x)=\sum_{k=0}^{N-1}x^{a_k},\ a_0=M>a_1>\dots>a_{N-1}=0,$$
let $\{\theta_1,\dots,\theta_{N-1}\}\subset(0,1)$ be a spectrum for $A(x)$, 
$$\epsilon_{jk}=e(\theta_j-\theta_k),\ \theta_0=0,$$
and let $z_1,\dots,z_M$ be the roots of the polynomial $A(x)$. The matrix
$(e(\theta_ia_j))_{i,j=0}^{N-1}$ is orthogonal. Therefore,
for $j\neq k$
$$\sum_{i=0}^{N-1}e(\theta_i(a_j-a_k))=0,$$
or
\begin{equation}\label{irr.e1}
  \sum_{i=1}^{N-1}e(\theta_i(a_j-a_k))=-1.
\end{equation}

Denote 
$$S_j=\sum_{i=1}^Mz_i^j.$$
Let $G$ be the Galois group of $A(x)$. Then,  by (\ref{irr.e1}),
for any $\sigma\in G$
$$\sum_{i=1}^{N-1}\sigma(e(\theta_i))^{a_j-a_k}=-1.$$
Averaging over $\sigma$, we get
\begin{equation}\label{irr.e2}
  S_{a_j-a_k}=-M/(N-1).
\end{equation}

By Newton's identities, if
$$A(x)=\sum_{j=0}^Mb_jx^j,$$
then 
\begin{equation}\label{irr.e3}
  S_j+b_{M-1}S_{j-1}+\dots+b_{M-j+1}S_1+jb_{M-j}=0.
\end{equation}
Taking consequently $j=1,\dots,M-a_1-1$, and using that all 
coefficients $b_i$ in (\ref{irr.e3}) are zeros, we get
\begin{equation}\label{irr.e4}
  S_1=\dots=S_{M-a_1-1}=0.
\end{equation}
Furthermore, for $j=M-a_1$ Newton's identity gives
\begin{equation}\label{irr.e5}
  S_{M-a_1}+(M-a_1)=0.
\end{equation}
On the other hand, $S_{M-a_1}=-M/(N-1)$ by (\ref{irr.e2}). Therefore,
\begin{equation}\label{irr.e6}
  M-a_1=M/(N-1).
\end{equation}

We claim that 
\begin{equation}\label{irr.e7}
  a_{j-1}-a_j\ge M/(N-1),\ j=1,\dots,N-1.
\end{equation}
Indeed, suppose the contrary. Then, by (\ref{irr.e4}), $S_{a_{j-1}-a_j}=0$,
but this equality does not agree with (\ref{irr.e2}). Hence,
$$M=\sum_{j=1}^{N-1}(a_{j-1}-a_j)\ge\sum_{j=1}^{N-1}M/(N-1)=M.$$
Thus, the inequalities in (\ref{irr.e7}) are actually equalities,
and we have
$$a_j=M-jM/(N-1),\ j=1,\dots,N-1,$$
and
$$A(x)=\sum_{j=0}^{N-1}x^{jM/(N-1)}
=\frac{x^{MN/(N-1)}-1}{x^{M/(N-1)}-1}.$$
In particular, all roots of $A(x)$ are roots of unity.  Since $A(x)$ is 
irreducible, $A(x)=\Phi_s(x)$ for some $s\in\nn$; moreover, $A(1)=N>1$
implies that $N=p$ and $s=p^\alpha$ for some prime $p$.  Hence 
$A(x)=(x^{p^\alpha}-1)/((x^{p^{\alpha-1}}-1)$ and
\[
A=\{0,p^{\alpha -1}, 2p^{\alpha -1},\dots,(p-1)p^{\alpha -1}\}.
\]
It is easy to see that $A$ tiles $\zz$ with the translation set
$B=\{0,1,\dots,p^{\alpha-1}-1\}+p^\alpha\zz$.


\section{Proof of Theorem \ref{N-factors}}
\label{sec-cube}
\init


We will consider polynomials of the form
\begin{equation}\label{cube-A}
A(x)=\prod_{i=1}^N A_i(x),\ 
A_i(x)=1+x^{m_i}+\dots+x^{m_i(n_i-1)}=\frac{x^{m_in_i}-1}{x^{m_i}-1}.
\end{equation}
It suffices to prove Theorem \ref{N-factors} under the assumption that
\begin{equation}\label{cube-m}
  (m_1,\dots,m_N)=1.
\end{equation}
Indeed, suppose that $(m_1,\dots,m_N)=d>1$, and let $A'=A/d$,
$m'_i=m_i/d$.  Then $A'$ has the form (\ref{cube-A}) and satisfies
(\ref{cube-m}).  Furthermore, $A$ tiles $\zz$ if and only $A'$
tiles $\zz$, and $A$ satisfies (T1)-(T2) if and only if so does $A'$
(see \cite{CM}).

Assume for now that $m_i,n_i$ are chosen so that $A(x)$ has $0,1$
coefficients.  (By Theorem \ref{N-factors},
(\ref{e-tower}) below is a sufficient condition.)

Let $m=(m_1,\dots,m_N)\in\rr^N$.
Consider the projection $\pi:\rr^N\to\rr$ given by
\[
\pi:\ (u_1,\dots,u_N)=u\to \jap{u,m}=u_1m_1+\dots+u_Nm_N.
\]
Let
\[
\cala=\{(j_1,\dots,j_N):\ j_k=0,1,\dots,n_k-1\}
\]
so that $\pi(\cala)=A$, and
\[
W=\{(w_1,\dots,w_N):\ w_i\in\zz,\ \jap{w,m}=0\}.
\]
If $A$ tiles $\zz$ with the translation set $B$, we will write
\[
\calb=\{(u_1,\dots,u_N):\ \jap{u,m}\in B\}=\pi^{-1}(B).
\]
Finally, we will denote $d_{ij}=(m_i,m_j)$.  We will sometimes
identify $\cala$ with the rectangular box $\{x\in\rr^N:\ 0\leq x_j
<n_j\}$.

\begin{lemma}\label{cube-lemma1}
Assume that $A\oplus B=\zz$, then:

(i) $\cala\oplus \calb$ is a tiling of $\zz^N$;

(ii) $\calb$ is invariant under all translations by vectors in $W$.
\end{lemma}

\proof
Let $w\in \zz^N$, then $\pi(w)=a+b$ for unique $a\in A, b\in B$.  
Let $u=\pi^{-1}(a)$; we are assuming that $\pi$ is one-to-one on
$\cala$, hence $u$ is uniquely determined.  Let also $v=w-u$.
Then $\pi(v)=\pi(w)-\pi(u)=b$, hence $v\in\calb$.  This shows that
each $w$ can be represented as $u+v$ with $u\in\cala,v\in\calb$.
Furthermore, for any such representation we must have $\pi(u)=a$ 
and $\pi(v)=b$, so that the above argument also shows uniqueness.
\qed

\bigskip

{\bf Remark} We also have the following converse of Lemma
\ref{cube-lemma1}.  Let a tiling $\cala\oplus\calb=\zz^N$ be given, where 
$\cala$ and $\calb$ are as above.  We claim that if (ii) holds, then
$A\oplus B=\zz$, where $A=\pi(\cala)$ and $B=\pi(\calb)$.
Indeed, by (\ref{cube-m}) $\pi$ is onto.  Let $x\in\zz$ and
pick a vector in $\pi^{-1}(x)$; this vector can be written as
$u+v$, where $u\in\cala$ and $v\in\calb$.  Therefore
$x=a+b$ with $a=\pi(u)\in A$ and $b=\pi(v)\in B$.  It remains 
to verify that this representation is unique.  Indeed, suppose that
$x=\pi(w)=\pi(w')$, then $\pi(w-w')=0$ so that $w-w'\in W$.  By
(ii), the tiling $\cala\oplus\calb$ is invariant under the translation
by $w-w'$.  Hence if we write $w=u+v$, $w'=u'+v'$ with 
$u,u'\in \cala$, $v,v'\in\calb$, it follows that $u=u'$ and consequently
$a=\pi(u)=\pi(u')$ is uniquely determined.  This also determines
$b=x-a$.

\bigskip

\begin{lemma}\label{span-lemma}
Let $w_{ij}$ be the vector whose $i$-th coordinate is $m_j/d_{ij}$,
$j$-th coordinate is $-m_i/d_{ij}$, and all other coordinates are 0.
Then
\[
W=\{\sum_{i,j}k_{ij}w_{ij}:\ k_{ij}\in\zz\}.
\]
\end{lemma}

\proof
Denote the set on the right by $W'$.  Since $w_{ij}$ have integer 
coordinates and $\jap{w_{ij},m}=0$, it is clear that $W'\subset W$.
We will now prove the converse using induction in $N$.  The inductive
step will not necessarily preserve the property (\ref{cube-m}).
However, if the lemma is proved for some $N$ under the assumption
(\ref{cube-m}), it also holds for the same $N$ without this assumption.
Indeed, suppose that $(m_1,\dots,m_N)=d>1$, then $d$ divides each $d_{ij}$,
so that we may replace each $m_j$ by $m'_j=m_j/d$ and apply the
version of the lemma in which (\ref{cube-m}) is assumed.  

The case $N=1$ is trivial since $\jap{w,m}=0$ in dimension 1 only if $w=0$.
Suppose that the lemma has been proved for $N-1$. 
We will show that any $w\in W$
can be written as $w=w'+w''$, where $w'\in W'$ and $w''\in W$, $w''_1=0$;
then the claim will follow by induction.
It suffices to prove that 
\[
w_1=\sum_{j=2}^N k_j\frac{m_j}{d_{1j}}
\]
for some choice of integers $k_j$; in other words, that
$(\frac{m_2}{d_{12}},\dots,\frac{m_N}{d_{1N}})$ divides $w_1$.
Since $\jap{w,m}=0$, we have
\[
m_1w_1=-m_2w_2-\dots -m_Nw_N.
\]
Hence $(m_2,\dots,m_N)$ divides $m_1w_1$.  By (\ref{cube-m}), it must in fact
divide $w_1$.  It only remains to observe that $(\frac{m_2}{d_{12}},\dots,
\frac{m_N}{d_{1N}})$ divides $(m_2,\dots,m_N)$.
\qed

\begin{theorem}\label{Nfactors-tiling}
Let $A$ be as in (\ref{cube-A}).  Then the following are equivalent:

\smallskip
(i) $A$ tiles $\zz$ by translations;

\smallskip
(ii) $A$ satisfies (T1)--(T2);

\smallskip
(iii) there is a labelling of the factors $A_i$ for which the following holds:
\begin{equation}\label{e-tower}
  \begin{array}{l}
  n_1|(\frac{m_2}{d_{12}}, \dots,\frac{m_N}{d_{1N}}),\\[3mm]
  n_2|(\frac{m_3}{d_{23}}, \dots,\frac{m_N}{d_{2N}}),\\[3mm]
  \dots,\\[3mm]
  n_{N-1}|\frac{m_N}{d_{N-1,N}}.  
  \end{array}
\end{equation}
\end{theorem}

Recall that we are assuming (\ref{cube-m}) throughout this section, 
including the proof that follows; however, it is easy to see that
the theorem remains true without this assumption (see the remark after
(\ref{cube-m})).

\medskip\noindent
{\it Proof of Theorem \ref{Nfactors-tiling}.}
We will prove that (i) $\Rightarrow$ (iii) $\Rightarrow$ (ii);
the implication (ii) $\Rightarrow$ (i) is proved (for more general $A$)
in \cite{CM}.

\medskip

\noindent {\bf (i) implies (iii):}
We will say that a set $V\subset \rr^N$ has {\em Keller's property} if
for each $v\in V$, $v\neq 0$, we have $v_i\in \zz\setminus\{0\}$
for at least one $i$. Let $L$ be the linear transformation on $\rr^N$ 
defined by 
\[
L(u_1,\dots,u_N)=(\frac{u_1}{n_1},\dots,\frac{u_N}{n_N}).
\]
If we identify $\cala$ with the rectangular box $\{x\in\rr^N:\ 0\leq x_j
<n_j\}$, then $L(\cala)$ is the unit cube $Q$ in $\rr^N$, and by Lemma
\ref{cube-lemma1}(i) $Q\oplus L(\calb)$ is a tiling of $\rr^N$.  We now
use the following theorem of Keller on cube tilings \cite{Keller}.

\begin{theorem}\label{thm-Keller}
\cite{Keller} If $Q\oplus V$ is a tiling of $\rr^N$, then the set 
$V-V:=\{v-v':\ v,v'\in V\}$ has Keller's property.
\end{theorem}

It follows that $L(\calb)-L(\calb)$ has Keller's property; in particular,
since $W\subset \calb-\calb$, $L(W)$ has Keller's property.  

We first claim that Keller's property for $W$ implies
the first equation in (\ref{e-tower}) for some labelling of $A_i$.  
Indeed, suppose that the first equation in (\ref{e-tower}) fails for
all such labellings.  Then for each $i\in\{1,\dots,N\}$ there is a 
$\sigma(i)\neq i$ such that $n_i\ndiv \frac{m_i}{d_{i\sigma(i)}}$.
We may find a cycle $i_1,\dots,i_r$ such that $i_{j+1}=\sigma(i_{j})$,
with $i_{r+1}=i_1$.  We thus have 
\begin{equation}\label{g.e1}
  n_{i_j}\ndiv \frac{m_{i_{j+1}}}{d_{i_{j},i_{j+1}}}
\end{equation}
for $j=1,\dots,r$.  

Define $w_{ij}$ as in Lemma \ref{span-lemma}.
If there is a $j$ such that
\begin{equation}\label{g.e2}
  n_{i_{j+1}}\ndiv \frac{m_{i_j}}{d_{i_{j},i_{j+1}}},
\end{equation}
then by (\ref{g.e1}), (\ref{g.e2}) Keller's property fails for $w_{i_{j},i_{j+1}}$.
If on the other hand (\ref{g.e2}) fails for all $j$, then this together
with (\ref{g.e1}) implies that Keller's property fails for 
$\sum_{j=1}^r w_{i_{j},i_{j+1}}$.  This completes the proof of the claim.

The remaining equations in (\ref{e-tower}) can now be obtained by induction
in $N$. Indeed, consider the set 
\[
W_1=\{(w_2,\dots,w_N): \ (0,w_2,\dots,w_N)\in W\}\subset\rr^{N-1}.
\]
This set (as a subset of $\rr^{N-1}$) has Keller's property, hence the previous
argument with $W$ replaced by $W_1$ implies the second equation in (\ref{e-tower}).
Similarly we obtain the rest of (\ref{e-tower}).

\medskip

\noindent {\bf (iii) implies (ii):}
By the definition of $A_i(x)$, 
\begin{equation}\label{cube-div}
 \Phi_s(x)\div A_i(x)\hbox{ if and only if }
 s\div m_in_i,\ s\ndiv m_i.
\end{equation}

We first prove (T1).  By the definition of $A_i(x)$, all its irreducible factors
are distinct cyclotomic polynomials, so that by (\ref{e.cyclo}) (T1) holds for 
each $A_i(x)$.  It therefore suffices to prove that if (\ref{e-tower}) holds,
then any prime power cyclotomic polynomial can divide at most one $A_i(x)$.
 
Let $p$ be a prime such that $\Phi_{p^\alpha}(x)$ divides 
$A_i(x)$ for some $\alpha,i$; it suffices to prove that $\Phi_{p^\alpha}(x)$ cannot
divide $A_j(x)$ for any $j>i$.  Let $p^{\beta_k}||m_k$ and $p^{\gamma_k}||n_k$
for $k=1,\dots,N$, then 
\begin{equation}\label{e-beta}
 \Phi_{p^\alpha}(x)\div A_k(x)\hbox{ if and only if }
 \beta_k<\alpha\leq\beta_k+\gamma_k.
\end{equation} 
In particular, it follows that $\gamma_i\neq 0$. 

Let $j>i$.  By (\ref{e-tower}) we have $n_i\div \frac{m_j}{d_{ij}}$, i.e.
$n_i(m_i,m_j)\div m_j$.  Thus $\gamma_i+\min(\beta_i,\beta_j)\leq \beta_j$.  
Note that we cannot have $\min(\beta_i,\beta_j)=\beta_j$, since then $\gamma_i$
would be 0. Hence $\min(\beta_i,\beta_j)=\beta_i$ and $\alpha\leq \beta_i+\gamma_i
\leq \beta_j$.  This and (\ref{e-beta}) imply that $\Phi_{p^\alpha}(x)\ndiv
A_j(x)$, as claimed.

We note for future reference that we have also proved the following:
\begin{equation}\label{e-tower2}
\hbox{if } \Phi_{p^\alpha}(x)\div A_i(x)\hbox{ for some }\alpha, 
\hbox{ then }\beta_i+\gamma_i\leq\beta_j\hbox{ for all }j>i.
\end{equation}

It remains to prove (T2).  We must prove that if $s>1$ is an integer such
that $\Phi_{p^\alpha}(x)\div A(x)$ for every $p^\alpha || s$, then
$\Phi_s(x)\div A(x)$.  We will in fact show that $\Phi_s(x)\div A_j(x)$,
where 
\[
j=\max\{k:\ \Phi_{p^\alpha}(x)\div A_k(x)\hbox{ for some }p^\alpha||s\}.
\]
By (\ref{cube-div}), it suffices to prove that $s\div m_jn_j$ and $s\ndiv m_j$.

For every $p^\alpha||s$ we have $\Phi_{p^\alpha}(x)\div A_k(x)$ for some $k\leq
j$.  Therefore $p^\alpha\div m_jn_j$; this follows from (\ref{cube-div}) if
$k=j$, and from (\ref{e-tower2}) if $k<j$.  Hence $s\div m_jn_j$.  On the other 
hand, by the definition of $j$ there is at least one prime power $p^\alpha||s$
such that $\Phi_{p^\alpha}(x)\div A_j(x)$.  By (\ref{cube-div}) we have
$p^\alpha\ndiv m_j$, so that $s\ndiv m_j$.
\qed

\section{Proof of Theorem \ref{two-factors}}
\label{sec3}
\init

In this section we will assume that $A(x)$ is as in (\ref{A1A2}).
Denote also $d=(m_1,m_2)$. We will prove that, under the above 
hypotheses, each of (T), (S), (T1)-(T2) is equivalent to the statement that
one of the following holds:
\begin{equation}\label{e.2main}
n_1\div\frac{m_2}{d},
\end{equation}
\begin{equation}\label{e.2mainb}
n_2\div\frac{m_1}{d}.
\end{equation}
We record for future reference that $\Phi_s(x)\div A(x)$ if and only if
\begin{equation}\label{cyclo-divide}
  s\div m_in_i,\ s\ndiv m_i\hbox{ for at least one of }i=1,2.
\end{equation}

By Theorem \ref{Nfactors-tiling} and the remark following it,
the statement that one of (\ref{e.2main}), 
(\ref{e.2mainb}) holds is equivalent to (T) and to (T1)-(T2).
In light of \cite{L}, Theorem 1.5,  this also implies (S).
It remains to show that 
(S) implies one of (\ref{e.2main}), (\ref{e.2mainb}).

Suppose that $A(x)$ has an $N$-spectrum
$\{\theta_j:\ j=1,\dots,N-1\}$.  
Let $\theta_N=0$, $\epsilon_j=e^{2\pi i\theta_j}\ (j=1,\dots,N-1)$,
$\epsilon_N=1$. Then the numbers
$$
\epsilon_j/\epsilon_k=e^{2\pi i(\theta_j-\theta_k)}
$$
are roots of $A(x)$ for all $j\neq k$, $j\le N$, $k\le N$.

We will first prove that one of the following must hold:
\begin{equation}\label{2.e11}
  \forall j\  m_1n_1\theta_j\in\zz,
\end{equation}
\begin{equation}\label{2.e12}
  \forall j\ m_2n_2\theta_j\in\zz.
\end{equation}
Indeed, suppose that (\ref{2.e11}) and (\ref{2.e12}) fail. Then there 
exist $j$ and $k$ such that
\begin{equation}\label{2.e13}
  m_1n_1\theta_j\not\in\zz,
\end{equation}
\begin{equation}\label{2.e14}
  m_2n_2\theta_k\not\in\zz.
\end{equation}
Since $\epsilon_j$ is a root of $A(x)$, we get from (\ref{2.e13}) that
\begin{equation}\label{2.e15}
 m_2n_2\theta_j\in\zz.
\end{equation}
Similarly,
\begin{equation}\label{2.e16}
  m_1n_1\theta_k\in\zz.
\end{equation}
The conditions (\ref{2.e13})--(\ref{2.e16}) imply
$$m_1n_1(\theta_j-\theta_k)\not\in\zz,$$
$$m_2n_2(\theta_j-\theta_k)\not\in\zz.$$
Thus, $\epsilon_j/\epsilon_k$ is not a root of $A(x)$. This contradiction 
shows that our supposition cannot occur. Without loss of generality,
we will assume that (\ref{2.e11}) holds.

For $l=0,\dots,n_1-1$ denote
$$J_l=\{j:m_1n_1\theta_j\equiv l\ (\bmod n_1)\}.$$
For $j,k\in J_l$, $j\neq k$, the number $\epsilon_j/\epsilon_k$
is not a root of $A_1(x)$. Hence, it is a root of $A_2(x)$. This means
that, for $j,k\in J_l$, $j\neq k$, the numbers $m_2n_2(\theta_j-\theta_k)$
are integers not divisible by $n_2$. This yields $|J_l|\le n_2$. On the 
other hand, the equality
$$N=n_1n_2=\sum_{l=0}^{n_1-1}|J_l|$$
demonstrates that actually $|J_l|= n_2$ for all $l$, and, moreover, for
a fixed $k\in J_l$, the numbers $m_2n_2(\theta_j-\theta_k)$ run over the 
complete residue system modulo $n_2$.

In particular, there exists $j\in J_0$ such that 
$$m_2n_2\theta_j\equiv1\ (\bmod n_2).$$ 
Therefore,
\begin{equation}\label{2.e17}
  {m_1m_2n_2 \over s}\theta_j\equiv{m_1\over s}\ (\bmod n_2).
\end{equation}
On the other hand, the condition $j\in J_0$ means $m_1\theta_j\in\zz$.
Therefore, 
$$m_1n_2\theta_j\equiv0\ (\bmod n_2)$$
and

\begin{equation}\label{2.e18}
  {m_1m_2n_2 \over s}\theta_j\equiv0\ (\bmod n_2).
\end{equation}
Comparing (\ref{2.e17}) and (\ref{2.e18}), we obtain (\ref{e.2mainb}).

\vskip.5in

\noindent{\bf Acknowledgement.} The first author was supported by Grants
02-01-00248 and 00-15-96109 from the Russian Foundation for Basic Research.
The second author was supported by the NSERC Grant 22R80520.



\bigskip
\noindent{\sc Department of Mechanics and Mathematics, Moscow State University,
Moscow 11992, Russia}

\smallskip
\noindent E-mail address: {\tt kon@mech.math.msu.su}

\bigskip
\noindent {\sc Department of Mathematics, University of British Columbia, Vancouver,
B.C. V6T 1Z2, Canada}

\smallskip
\noindent E-mail address: {\tt ilaba@math.ubc.ca}

\end{document}